\documentclass[11pt,leqno,twoside]{article}

\usepackage{amsfonts,amsmath,amsthm,amssymb}
\usepackage{enumerate}
\usepackage{graphics,graphicx,subfigure}

\usepackage{color}
\usepackage{todonotes}
\usepackage{cancel}
\usepackage{url}
\usepackage{hyperref}
\usepackage{makeidx}
\usepackage{showidx}
\usepackage{multicol}        
\usepackage{xspace}
\usepackage{stmaryrd}        
\usepackage{pifont}          
\usepackage{fancybox}        
\usepackage{bm}

 \usepackage{fancyhdr}
\theoremstyle{plain}
 \theoremstyle{definition}
 \newtheorem{lem}{Lemma}
 \newtheorem{defn}[lem]{Definition}
 \newtheorem{thm}[lem]{Theorem}
 \newtheorem{prop}[lem]{Proposition}
 \newtheorem{cor}[lem]{Corollary}
 \newtheorem{notn}[lem]{Notations}
 \newtheorem{pb}[lem]{Problem}
 \newtheorem{form}[lem]{Formulae}
 
 \newtheorem*{rk}{Remark}
 \newtheorem*{com}{Comment}
 \newtheorem*{ex}{Example}
 \theoremstyle{remark}

 \newcommand{\blem}{\begin{lem}}
 \newcommand{\elem}{\end{lem}}
 \newcommand{\bdefn}{\begin{defn}}
 \newcommand{\edefn}{\end{defn}}
 \newcommand{\bthm}{\begin{thm} }
 \newcommand{\ethm}{\end{thm}}
 \newcommand{\bprop}{\begin{prop}}
 \newcommand{\eprop}{\end{prop}}
 \newcommand{\bcor}{\begin{cor}}
 \newcommand{\ecor}{\end{cor}}
 \newcommand{\bnotn}{\begin{notn}}
 \newcommand{\enotn}{\end{notn}}
 \newcommand{\bpb}{\begin{pb}}
 \newcommand{\epb}{\end{pb}}
 \newcommand{\bform}{\begin{form}}
 \newcommand{\eform}{\end{form}}
 \newcommand{\brk}{\begin{rk}}
 \newcommand{\erk}{\end{rk}}
 \newcommand{\bcom}{\begin{com}}
 \newcommand{\ecom}{\end{com}}
 \newcommand{\bex}{\begin{ex}}
 \newcommand{\eex}{\end{ex}}
 \newcommand{\bpf}{\begin{proof}}
 \newcommand{\epf}{\end{proof}}

\newcommand{\cC}{\mathcal{C}}

\newcommand{\cK}{\mathcal{K}}

\newcommand{\cP}{\mathcal{P}}

\newcommand{\cV}{\mathcal{V}}

\newcommand{\cX}{\mathcal{X}}

\newcommand{\bR}{\mathbb{R}}

\newcommand{\bT}{\mathbb{T}}

\newcommand{\be}{\begin{equation}}
\newcommand{\ee}{\end{equation}}
\newcommand{\bal}{\begin{align}}
\newcommand{\eal}{\end{align}}
\newcommand{\ba}{\begin{align*}}
\newcommand{\ea}{\end{align*}}
\newcommand{\bmx}{\begin{matrix}}
\newcommand{\emx}{\end{matrix}}
\newcommand{\bbmx}{\begin{bmatrix}}
\newcommand{\ebmx}{\end{bmatrix}}
\newcommand{\bpmx}{\begin{pmatrix}}
\newcommand{\epmx}{\end{pmatrix}}
\newcommand{\bvmx}{\begin{vmatrix}}
\newcommand{\evmx}{\end{vmatrix}}

\newcommand{\ul}{\underline}
\newcommand{\ol}{\overline}
\newcommand{\wh}{\widehat}
\newcommand{\wt}{\widetilde}
\newcommand{\f}{\frac}

\newcommand{\inc}{\subseteq}

\newcommand{\setm}{\setminus}

\newcommand{\Id}{\mathrm{Id}}

\newcommand{\sgn}{\mathrm{sgn}}

\newcommand{\minimize}[1]{\underset{#1}{\rm minimize}\,}

\newcommand{\la}{\lambda}
\newcommand{\La}{\Lambda}
\newcommand{\eps}{\varepsilon}
  
\newcommand{\rev}[1]{{#1}}

\usepackage{framed}
\usepackage{algorithm}
\usepackage{algpseudocode}

\title{\vspace{-15mm}
Full Recovery from Point Values: \\an Optimal Algorithm for Chebyshev Approximability Prior
\medskip\hrule height 1.2pt \vspace{-6mm}}
\author{Simon Foucart\footnote{S. F. partially supported by grants from the NSF (DMS-2053172) and from the ONR (N00014-20-1-2787).
\url{foucart@tamu.edu} }  --- Texas A\&M University}
\date{\vspace{-6mm}\rule{100mm}{0.8pt}}

\newcommand\shorttitle{Full Recovery from Point Values}
\newcommand\authors{S. Foucart}

\begin{document}
\maketitle

\vspace{-15mm}
\begin{abstract}\vspace{-3mm}
Given pointwise samples of an unknown function belonging to a certain model set,
one seeks in Optimal Recovery to recover this function in a way that minimizes the worst-case error of the recovery procedure.
While it is often known that such an optimal recovery procedure can be chosen to be linear,
e.g. when the model set is based on approximability by a subspace of continuous functions,
a construction of the procedure is rarely available.
This note uncovers a practical algorithm to construct a linear optimal recovery map
when the approximation space is a Chevyshev space \rev{of univariate functions that has dimension at least three and contains the constants}. 
\end{abstract}

\noindent {\it Key words and phrases:}  Optimal recovery, Chebyshev spaces,  $\ell_1$-minimization, simplex algorithm.

\noindent {\it AMS classification:} 41A05, 41A10, 41A50, 90C05.

\vspace{-5mm}
\begin{center}
\rule{100mm}{0.8pt}
\end{center}


\section{Problem Setting}

Throughout this note, one works in the space $C(\cX)$ of continuous functions on a compact set $\cX$ 
equipped with the uniform norm defined for $f \in C(\cX)$ by $\|f\|_{C(\cX)} = \max\{ |f(x)|, x \in \cX \}$.
Given points $x^{(1)},\ldots,x^{(m)} \in \cX$,
an unknown function $f \in C(\cX)$ is observed via the point values
$$
y_i = f(x^{(i)}),
\qquad i \in  \{ 1, \ldots,m \}.
$$
This so-called {\em a posteriori} information alone is not enough to approximate/learn/recover $f$ in any meaningful way.
One also needs some {\em a priori} information,
usually expressed by the membership of $f$ to some model set $\cK$,
i.e., by $f \in \cK$.
The performance of a recovery procedure---which is nothing else than a map $\Delta$ from $\bR^m$ into $C(\cX)$---can then be assessed via its worst-case error over~$\cK$, defined as
\be
\label{GWCE}
{\rm wce}_{\cK}(\Delta) := 
\sup_{f \in \cK} \big\|f - \Delta( [f(x^{(1)}); \ldots;  f(x^{(m)})] ) \big\|_{C(\cX)}.
\ee
\newpage
The question being addressed in this note reads as follows:
\begin{framed}\noindent \vspace{-5mm}
\begin{equation}
\tag{{\bf Q}}
\label{Q}
\mbox{Can one construct an {\rm optimal} recovery procedure $\Delta$,
i.e.,  one that minimizes ${\rm wce}_{\cK}(\Delta)$?}\vspace{-3mm}
\end{equation}
\end{framed}
This objective is too ambitious for a general model set $\cK$,
so one concentrates in particular on model sets based on approximation capabilities.
Precisely, given a subset $\cV$ of $C(\cX)$ and a parameter $\eps \ge 0$,
one considers the approximability model
$$
\cK_{\cV,\eps} = \big\{ f \in C(\cX): {\rm dist}_{C(\cX)}(f,\cV) \le \eps \big\}.
$$
The premise that the observed function $f$ belongs to $\cK_{\cV,\eps}$ can be viewed as 
making explicit an assumption often appearing implicitly in numerical methods.
Indeed,  with $\cP_n$ denotes the space of polynomials of degree $<n$,
one often targets error bounds featuring  ${\rm dist}_{C(\cX)}(f,\cP_n)$---hence no error when $f \in \cP_n$: this is the exactness principle---and this presupposes that ${\rm dist}_{C(\cX)}(f,\cP_n)$ should be small.
A typical example is supplied by the design of quadrature formulas, discussed in \cite{Tre} along the lines of the exactness principle.
In the rest of this note,
the implicit-made-explicit assumption therefore takes the form of the prior ${\rm dist}_{C(\cX)}(f,\cV) \le \eps$
for some $n$-dimensional subspace $\cV$ of $C(\cX)$ sharing key similitudes with $\cP_n$\rev{, i.e., }Chebyshev spaces containing constant functions.
Some important properties of Chebyshev spaces are recalled in Section \ref{SecChe}.
For now, one only mentions that these spaces essentially do not exist in dimension $d>1$.
Arguably, this restricts the impact of the positive answer to Question \eqref{Q}
given for the case $\cX = [-1,1]$.
However,  even in this elementary case,
a knowledge gap is still filled by the complete answer exposed in this note.
For sure, some pieces were known---they are recalled in Section \ref{SecOR}---but they did not provide a genuinely constructive recovery procedure.
Here,  a practical recovery algorithm is indeed provided.
The correctness of this algorithm is justified in Section \ref{SecJusti}.
Section \ref{SecRks} concludes with some related remarks,
including a recipe to compute the maximum ratio of uniform and discrete norms over Chebyshev spaces.

\section{Description of an Optimal Algorithm}
\label{SecAlgo}

Before launching into theoretical considerations,
one directly puts forward the proposed procedure as Algorithm \ref{Algo} below,
with points requiring justification indicated by a triangle $\triangleright$.
But first, one quickly elucidates the notation $a_S$ and $M_S$
for a vector $a \in \bR^m$ and for a matrix $M \in \bR^{n \times m}$
when $S \inc \{1,\ldots,m\}$ is an index set of size $|S| = n$:
they represent the subvector in $\bR^n$, resp. the submatrix  in $\bR^{n \times n}$,
obtained by keeping only entries, resp.  columns,  indexed by $S$.
Note that $n \le m$ is assumed, 
for otherwise there would exist $h \in \cV\setm \{0\} $ such that $h(x^{(1)}) = \cdots = h(x^{(1)})  = 0$,
which, by considering $f+th$ as $t \to \infty$,  would yield ${\rm wce}_{\cK_{\cV,\eps}}(\Delta)= \infty$ for any $\Delta$.

\begin{algorithm}[h]
\caption{Optimal recovery procedure from point values for Chebyshev approximability prior}\label{Algo}
\begin{algorithmic}
\Require Points $x^{(1)}, \ldots,  x^{(m)} \in [-1,1]$ (completed with $x^{(0)}=-1$ and $x^{(m+1)}=1$ if necessary)
and functions $v_1,\ldots,v_n$ forming a basis for the $n$-dimensional Chebyshev space $\cV \inc C[-1,1]$ containing the constant functions
\Ensure $m \ge n \ge 3$
\State
{\bf do} create the matrix $M \in \bR^{n \times m}$ having entries $M_{j,i} = v_j(x^{(i)})$,
$j \in \{1,\ldots,n \}$, $i \in \{ 1,\ldots,m \}$
\For{$k \in \{1,\ldots,m+1\}$ } 
\State 
select a point $z^{(k)} \in (x^{(k-1)},x^{(k)})$
\State
create the vector $b^{(k)} \in \bR^n$ having entries $b^{(k)}_j = v_j(z^{(k)})$, 
$j \in \{1,\ldots,n \}$
\State
\ul{compute} an index set $S_k \inc \{1,\ldots,m\}$ of size $|S_k | = n$ so that
the vector $a^{(k)} \in \bR^m$ supported on $S_k$ and satisfying $a^{(k)}_{S_k} =  M_{S_k} ^{-1} b^{(k)}$ 
is a solution to \Comment{}
$$
\minimize{a \in \bR^m} \; \|a\|_1 
\qquad \mbox{subject to } M a = b^{(k)}
$$
\EndFor
\For{$i \in \{1,\ldots,m\}$}
\State
define a \ul{continuous} function $a^\sharp_i$ on $[-1,1]$ by \Comment{}
\begin{align*}
a^\sharp_i(x) & = 
\left\{ \bmx \sum_{j=1}^n  \big[ M_{S_k}^{-1} \big]_{i,j} v_j(x)
& \mbox{ if } i \in S_k\\
0 & \mbox{ if } i \not\in S_k \emx \right\},
& x \in (x^{(k-1)},x^{(k)}), \qquad & k \in \{1,\ldots,m+1\}\\
a^\sharp_i(x^{(\ell)}) & = \left\{ \bmx 1 & \mbox{ if } i =\ell \\ 0 & \mbox{ if } i \not= \ell
\emx \right\}, & & \ell \in \{1,\ldots,m\}
\end{align*}
\EndFor
\State
\Return the \ul{optimal} recovery map $\Delta^\sharp: y \in \bR^m \mapsto \sum_{i=1}^m y_i a^\sharp_i \in C[-1,1]$
 \Comment{}
\end{algorithmic}
\end{algorithm}


To confirm that Algorithm \ref{Algo} acts as intended, 
three points need to be accounted for:
the existence and computability of an index set $S_k$ with the required property,
the continuity of the functions $a^\sharp_1,\ldots,a^\sharp_m$,
and the optimality of the recovery map $\Delta^\sharp$.
The last two points rely on results about Chebyshev spaces and Optimal Recovery,
which are covered in Sections \ref{SecChe} and \ref{SecOR}.
As for the first point,  it can be explained right now.
Concerning existence,  recall that an optimization program
\begin{equation} \label{L1Min}
\minimize{a \in \bR^m} \; \|a\|_1 
\qquad \mbox{subject to } M a = b
\end{equation}
always admits an $n$-sparse solution $\wh{a} \in \bR^m$
(this can be proved along the lines of \cite[Theorem~12.7]{Mal} and \cite[Theorem~6.1]{BookDS}),
say supported on some $S \inc \{1,\ldots,m\}$ with $|S|=n$.
The constraint $M \wh{a} = b$ then reads  $M_{S} \wh{a}_S = b $,
i.e., $\wh{a}_S = M_{S}^{-1} b$.
Concerning computability, it can be realized by recasting \eqref{L1Min} as a standard-form linear program  to be solved via the simplex algorithm.  
Precisely,  introducing a slack variable $c = [a^+;a^-] \in \bR^{2m}$ with nonnegative vectors $a^+,a^- \in \bR^{m}$ satisfying $a= a^+ - a^-$ and $|a| = a^+ + a^-$,
the $\ell_1$-minimization \eqref{L1Min}
is equivalent to
$$
\minimize{c \in \bR^{2m}} \; \sum_{j=1}^{2m} c_j 
\qquad \mbox{subject to } \big[ \; M \;  \big| -M \; \big] c = b
\; \mbox{ and } \; c \ge 0.
$$
Solving the latter with the simplex algorithm yields a solution $\wh{c} = [\wh{a}^+; \wh{a}^-]$ which is an extreme point of the feasibility polytope and as such (see e.g.  \cite[Lemma 20.2]{BookDS})
is $n$-sparse.
In turn, the solution $\wh{a} = \wh{a}^+ - \wh{a}^-$ to \eqref{L1Min} is also $n$-sparse. 

\section{Reminders on Chebyshev Spaces}
\label{SecChe}

Chebyshev subspaces are at the center of Approximation Theory:
they are defined as the subspaces from which best approximants always exist and are unique.
Remarkably,
in $C(\cX)$,
they coincide with the subspaces for which Lagrange interpolation is always possible and unique.
Precisely, an $n$-dimensional subspace $\cV$ of $C(\cX)$ is a Chebyshev space if and only if,
for any distinct points $\xi^{(1)},\ldots,\xi^{(n)} \in \cX$ and any values $\gamma_1,\ldots,\gamma_n \in \bR$,
there exists a unique function $v \in \cV$ such that $v(\xi^{(1)}) = \gamma_1, \ldots,v(\xi^{(n)}) = \gamma_n$.  
By considering the linear map $v \in \cV \mapsto [v(\xi^{(1)}); \ldots; v(\xi^{(n)})] \in \bR^n$,
this is easily seen to be equivalent to the invertibility of the $n \times n$ matrix with entries $v_j(\xi^{(i)})$,
where $(v_1,\ldots,v_n)$ denotes a basis for~$\cV$.
Fixing such a basis, the determinant of this matrix must be nonzero,
and hence must be either always positive or always negative for all pointsets ${\bm \Xi}= (\xi^{(1)},\ldots,\xi^{(n)})$ satisfying $\xi^{(1)} < \cdots < \xi^{(n)}$.
Protected against division by zero,
one can now consider the function $L_{{\bm \Xi},i} \in \cV$ defined for $x \in \cX$ by
$$
L_{{\bm \Xi},i}(x) \hspace{-1mm}=\hspace{-1mm} 
\bvmx
 \; \cdots & v_1(\xi^{(i-1)}) & v_1(x) & v_1(\xi^{(i+1)}) & \cdots \; \\
 & \vdots & \vdots & \vdots & \\
 \; \cdots & v_n(\xi^{(i-1)}) & v_n(x) & v_n(\xi^{(i+1)}) & \cdots \;
\evmx
\Bigg/
\bvmx
\; \cdots & v_1(\xi^{(i-1)}) & v_1(\xi^{(i)}) & v_1(\xi^{(i+1)}) & \cdots \; \\
 & \vdots & \vdots & \vdots & \\
\; \cdots & v_n(\xi^{(i-1)}) & v_n(\xi^{(i)}) & v_n(\xi^{(i+1)}) & \cdots \; 
\evmx.
$$
It is called the $i$th fundamental Lagrange interpolator on ${\bm \Xi}$,
by virtue of $L_{{\bm \Xi},i}(\Xi^{(j)}) = \delta_{i,j}$ for $i,j \in \{1,\ldots,n\}$.
Note that $L_{{\bm \Xi},i}$ has no zeros besides $\xi^{(1)}, \ldots, \xi^{(i-1)},\xi^{(i+1)},\ldots,\xi^{(n)}$,
otherwise $\gamma_1 = \cdots = \gamma_n = 0$ could be interpolated on a set of $n$ distinct points by two different functions from $\cV$, namely by $0$ and by $L_{{\bm \Xi},i}$.
Note also that, for a fixed $x \not\in \{ \xi^{(1)}, \ldots, \xi^{(n)} \}$,
the sequence  $( L_{{\bm \Xi},i}(x)  )_{i=1}^n$ cannot keep a constant sign when $n \ge 3$:
if $x \in (\xi^{(i)},\xi^{i+1})$ for some $i \in \{1,\ldots,n-1\}$,
then $L_{{\bm \Xi},i}(x)$ and $L_{{\bm \Xi},i+1}(x)$ are both positive,
but $L_{{\bm \Xi},i-1}(x)$ or $L_{{\bm \Xi},i+2}(x)$---whichever exists---is negative,
and if $x < \xi^{(1)}$, say, 
then $L_{{\bm \Xi},1}(x)$ is positive, but  $L_{{\bm \Xi},2}(x)$ is negative.

In the multivariate situation,  it is easy to realize that the matrix with entries $v_j(\xi^{(i)})$ cannot be invertible for all choices of distinct points $\xi^{(1)},\ldots,\xi^{(n)}$. 
This is formalized by Mairhuber--Curtis theorem:
if $ d >1 $ and $\cX \inc \bR^d$ contains an interior point, 
then there is no Chebyshev subspace $\cV$ of $C(\cX)$ with dimension $n \ge  2$.
Thus, one usually considers Chebyshev spaces on $\cX = [-1,1]$ and $\cX = \bT$,
with the prototypical examples being spaces of algebraic polynomials and  of trigonometric polynomials.
In these univariate situations,  involving a simple notion of differentiation,
Lagrange interpolation can painlessly be generalized to Hermite interpolation,
leading to the introduction of extended Chebyshev spaces as spaces for which Hermite interpolation is always possible and unique.
On compact intervals $\cX = [a,b]$,
 extended Chebyshev spaces even turn out to be extended complete Chebyshev spaces\footnote{This result can be found in \cite[Theorem 5, p 97]{Cop} with a different terminology. 
 It has also been reproved in a simpler way in \cite[Appendix A.3]{Thesis}.}.
 These spaces do not need a formal definition here, 
 as it suffices to say that they are characterized by the existence of a basis $(u_0,\ldots,u_{n-1})$ of the form 
\begin{align*}
u_0(x) &= w_0(x),\\
u_1(x) & = w_0(x) \int_t^x w_1(x_1) dx_1,\\
\vdots & \\
u_{n-1}(x) & = w_0(x) \int_t^x w_1(x_1) \int_t^{x_1} w_2(x_2) \cdots \int_t^{x_{n-2}} w_{n-1}(x_{n-1}) \, dx_{n-1} \cdots dx_2 dx_1,
\end{align*}
relative to positive weights $w_0 \in C^{n-1}[a,b],  w_1 \in C^{n-2}[a,b], \ldots,w_{n-1} \in C[a,b]$ and a point $t \in [a,b]$.
Clearly, taking $w_0=\cdots=w_{n-1}=1$ generates the shifted monomial  basis with $u_j(x) = (x-t)^j/j!$.

\section{Reminders on Optimal Recovery}
\label{SecOR}

The question addressed in this note is an instance of the generic Optimal Recovery problem,
on which a brief rundown is laid out here.
The following results,  either classical or due to \cite{DFPW} for the approximability model,
can all be found in \cite[Chapters 9 and 10]{BookDS}.
In an abstract setting,
an object $f$ from a normed space $F$---not necessarily a function space---is assumed to belong to a model set $\cK \inc F$
and is observed via $y = \La f$ for some linear map $\La: F \to \bR^m$,
i.e., via $y_i = \la_i(f)$, $i = 1,\ldots,m$,
for some linear functionals $\la_1,\ldots,\la_m \in F^*$.
The goal is to recover not necessarily $f$ itself,
but $Q(f)$ for some linear map $Q: F \to Z$,
and to do so in an  optimal way, so as to minimize the worst-case error  over $\cK$,  defined as
$$
{\rm wce}_{\cK,Q}(\Delta) \rev{:=} \sup_{f \in \cK} \big\| Q(f) - \Delta(\La f) \big\|_Z .
$$
If the model set $\cK$ is symmetric and convex
and if the quantity of interest $Q: F \to \bR$ is a linear functional,
it is well known that the set of recovery maps $\Delta: \bR^m \to \bR$ minimizing ${\rm wce}_{\cK,Q}(\Delta)$ contains a linear map $\Delta^\sharp$.
This linearity result is a typical statement in Optimal Recovery (see \cite{Pac} for an in-depth discussion),
but it is not always constructive.
In case of the approximability set~\rev{$\cK_{\cV,\eps} = \{ f \in F : {\rm dist}_F(f,\cV) \le \eps \}$ relative to a linear subspace $\cV$ of $F$ and a parameter $\eps \ge 0$},
such a linear optimal recovery map is given as $\Delta^\sharp = \langle a^\sharp,  \cdot \rangle$,
where $a^\sharp \in \bR^m$ is a solution to 
\be
\label{MinAbs}
\minimize{a \in \bR^m} \; \bigg\| Q - \sum_{i=1}^m a_i \la_i \bigg\|_{F^*} 
\qquad \mbox{subject to } \sum_{i=1}^m a_i {\la_i}(v) = Q(v) \; \mbox{ for all } v \in \cV.
\ee
This result was proved in \cite{DFPW} as a consequence of the Hahn--Banach extension theorem. 
It can also be deduced from the above linearity result---which is a consequence of the Hahn--Banach separation theorem---using the following argument already outlined in \cite{Ins}:
since it is enough to minimize the worst-case error over $\cK_{\cV,\eps}$
 among linear recovery maps $\Delta_a =\langle a, \cdot \rangle$,
 and since the worst-case error for linear maps $\Delta_a$ is
\begin{align}
\nonumber
{\rm wce}_{\cK_{\cV,\eps},Q}(\Delta_a)
& = \sup_{f \in C(\cX)} \bigg\{ \bigg| Q(f) - \sum_{i=1}^m a_i \la_i(f) \bigg| : \|f-v\|_{C(\cX)} \le \eps \mbox{ for some } v \in \cV \bigg\}\\
\nonumber
& =   \sup_{\substack{f \in C(\cX) \\ v \in \cV}} \bigg\{ \bigg| Q(f) - \sum_{i=1}^m a_i \la_i(f) \bigg| : \|f-v\|_{C(\cX)} \le \eps \bigg\}\\
\nonumber
& = \sup_{\substack{g \in C(\cX) \\ v \in \cV}} \bigg\{ \bigg| \bigg( Q(g) - \sum_{i=1}^m a_i \la_i(g) \bigg) +  \bigg( Q(v) - \sum_{i=1}^m a_i \la_i(v) \bigg)  \bigg| : \|g\|_{C(\cX)} \le \eps \bigg\}\\
\nonumber
& = \sup_{g \in C(\cX) } \bigg\{ \bigg|  Q(g) - \sum_{i=1}^m a_i \la_i(g)  \bigg| : \|g\|_{C(\cX)} \le \eps \bigg\}
+ \sup_{v \in \cV} \bigg\{ \bigg|  Q(v) - \sum_{i=1}^m a_i \la_i(v)  \bigg|  \bigg\}\\
\label{WCE_lin}
& = \bigg\|  Q - \sum_{i=1}^m a_i \la_i  \bigg\|_{F^*} \times \eps
+ \left\{ \bmx
+\infty & \mbox{if $Q(v) \not= \sum_{i=1}^m a_i \la_i(v)$ for some $v \in \cV$}\\
0 & \mbox{if $Q(v) = \sum_{i=1}^m a_i \la_i(v)$ for all $v \in \cV$ \; \;} 
\emx \right\}, 
\end{align}
the minimization of the latter among all $a \in \bR^m$ indeed reduces to the program \eqref{MinAbs}. 
 
It is also worth pointing out that the minimal worst-case error---aka intrinsic error---cannot exceed the so-called null error. 
Over the model set $\cK_{\cV,\eps}$,
this means that the minimal worst-case error for any linear quantity of interest $Q: F \to Z$ is lower-bounded by the product of the approximability parameter $\eps$ and an indicator $\mu_{\cV,Q}(\La)$ of the compatibility between the model (through $\cV$) and the observation process (through $\La$).
Precisely, one has 
\be
\label{MuEps}
\min_{\Delta: \bR^m \to \bR} {\rm wce}_{\cK_{\cV,\eps},Q}(\Delta) 
\ge  \mu_{\cV,Q}(\La) \times \eps,
\qquad \mbox{ where }
\mu_{\cV,Q}(\La) := \sup_{h \in \ker(\La) \setm\{0\}} \f{\|Q(h)\|_Z}{{\rm dist}_{C(\cX)}(h,\cV)},
\ee
with equality occurring when $Q: F \to \bR$ is a linear functional.
\rev{The indicator $\mu_{\cV,Q}(\La)$ was introduced in the article  \cite{DFPW},
which provides references to its earlier appearances in the case $Q=\Id$,
where it can be interpreted as the reciprocal of an angle between $\cV$ and $\ker(\La)$.} 

In the framework of this note,
where the space $F$ is $C(\cX)$ and
the observation functionals $\la_i $ are evaluations $\delta_{x^{(i)}}$ at points $x^{(i)}$,
if $Q = \delta_x$ is the evaluation at a point $x \not\in  \{x^{(1)}, \ldots, x^{(m)} \}$,
then one has $\| Q - \sum_{i=1}^m a_i \la_i \|_{F^*} = 1+ \|a\|_1$,
hence an optimal recovery map over $\cK_{\cV,\eps}$ takes the form $\Delta^\sharp = \langle a^\sharp(x) , \cdot \rangle$,
where $a^\sharp(x) \in \bR^m$ is a solution to 
\be
\label{MinSpe}
\minimize{a \in \bR^m} \; \|a\|_1
\qquad \mbox{subject to } \sum_{i=1}^m a_i v(x^{(i)}) = v(x) \; \mbox{ for all } v \in \cV.
\ee
\rev{Here and in the rest of this section,
the subspace $\cV \inc C(\cX)$ is arbitrary
and the following facts do not yet rely on it being a Chebyshev space.
For instance, 
the above optimization program always reduces to \eqref{L1Min} with $b = b(x):= [v_1(x);\ldots;v_n(x)]$,
simply by remarking that 
the constraint in \eqref{MinSpe} is met for all $v \in \cV$ if and only if it is met for all of the elements  of a basis $(v_1,\ldots,v_n)$ for~$\cV$.}
Furthermore,  since $Q=\delta_x$ is a linear functional,
\rev{ equality holds in \eqref{MuEps},
which,  in conjunction with~\eqref{WCE_lin}, 
leads to an important identity already found in \cite[Subsection 4.2]{DFPW}, namely}
$$
\mu_{\cV, \delta_x}(\La_{\bm X}) = 1 + \|a^\sharp (x) \|_1,
$$
where $\La_{\bm X}: F \to \bR^m$ denotes the linear map defined by $\La_{\bm X}(f) = [f(x^{(1)}); \ldots; f(x^{(m)})]$. 

For the full \rev{recovery} problem, i.e., for $Q = \Id$,
it was noticed in \cite{DFPW} that solving \eqref{MinSpe} for all $x \in \cX$ 
demonstrates the existence of an optimal recovery map $\Delta^\sharp: \bR^m \to C(\cX)$
which is linear.
Precisely,  if $a^\sharp(x) \in \bR^m$ represents again a (not necessarily unique) solution to \eqref{MinSpe},
then the linear map $\Delta^\sharp$ defined by $\Delta^\sharp(y)(x) = \sum_{i=1}^m y_i a^\sharp_i (x)$ minimizes the worst-case error over $\cK_{\cV,\eps}$.
Indeed,  for any $f \in \cK_{\cV,\eps}$, 
with $v \in \cV$ chosen so that ${\rm dist}_{C(\cX)}(f,\cV) = \|f-v\|_{C(\cX)}$,
one has
\begin{align*}
\|f - \Delta^\sharp (\La_{\bm X} f) \|_{C(\cX)}
& = \max_{x \in \cX} \Big| f(x) - \sum_{i=1}^m f(x^{(i)}) a^\sharp_i(x) \Big|
= \max_{x \in \cX} \Big| (f-v)(x) - \sum_{i=1}^m (f-v)(x^{(i)}) a^\sharp_i(x) \Big|\\
& \le \max_{x \in \cX}  \Big\| \delta_x - \sum_{i=1}^m a^\sharp_i(x) \delta_{x^{(i)}} \Big\|_{C(\cX)^*} \times \|f-v\|_{C(\cX)}
\le \max_{x \in \cX} \big(1+\|a^\sharp(x)\|_1 \big) \times \eps\\
& = \max_{x \in \cX} \mu_{\cV,\delta_x} (\La_{\bm X}) \times \eps. 
\end{align*}
Taking the supremum over  $f \in \cK_{\cV,\eps}$ while remarking $\max_{x \in \cX} \mu_{\cV,\delta_x} (\La_{\bm X}) = \mu_{\cV,\Id}(\La_{\bm X})$,
it follows with the help of \eqref{MuEps} that
$$
{\rm wce}_{\cK_{\cV,\eps}, \Id}(\Delta^\sharp)
 \le \mu_{\cV,\Id}(\La_{\bm X}) \times \eps
\le \min_{\Delta: \bR^m \to \bR} {\rm wce}_{\cK_{\cV,\eps},\Id}(\Delta).
$$
This establishes the optimality of $\Delta^\sharp$,
provided $\Delta^\sharp$ maps in $\cC(\cX)$,
i.e.,  provided $x \in \cX \mapsto a^\sharp(x) \in \bR^m$ can be made continuous by properly selecting the minimizers $a^\sharp(x)$ for all $x \in \cX$.
This was the intricate part of the argument in \cite{DFPW},
carried out under the proviso that the space $\cV$ contains the constant functions.
But evidently,  solving \eqref{MinSpe} for all $x \in \cX$ does not constitute a practical algorithm,
as opposed to Algorithm \ref{Algo}.
Still, the above considerations are the basis of the validation of Algorithm \ref{Algo}.
Note in passing that\rev{, for an arbitrary subspace $\cV$, it could be practical to solve \eqref{MinSpe} for all $x$ in a fine discretization of $\cX$.
This would not answer to Question \eqref{Q} satisfactorily,
but could nonetheless provide a good ersatz}.

\section{Validation of the Proposed Algorithm}
\label{SecJusti}

It is now time to validate Algorithm \ref{Algo} by showing that it does indeed return an optimal recovery map---a linear one, to boot---for the full recovery problem in $C[-1,1]$ over $\cK_{\cV,\eps}$ when $\cV$ is a Chebyshev space of dimension $n \ge 3$ and containing the constant functions.
Recall that proving this statement amounts to justifying the last two points indicated by some  $\triangleright$ in Algorithm \ref{Algo},
i.e., the continuity of $a^\sharp_1,\ldots,a^\sharp_m$
and the optimality of the recovery map $\Delta^\sharp$.
This task will rely on a characterization,  for any $x \in [-1,1]$,   of an $n$-sparse minimizer of the program \eqref{MinSpe}, written here  as
\be
\label{MinSpe2}
\minimize{a \in \bR^m} \|a\|_1 \qquad
\mbox{subject to } Ma = b(x),
\ee
recalling that the matrix $M \in \bR^{n \times m}$ and the vector $b(x) \in \bR^n$ have entries
$$
M_{j,i} = v_j(x^{(i)})
\qquad \mbox{and} \qquad 
b_j(x) = v_j(x).
$$
The characterization uses (in one direction only) the following simple observation.

\blem
\label{LemNZ}
Let $\cV$ be an $n$-dimensional Chebyshev subspace of $\cC[-1,1]$.
For $x \not\in \{x^{(1)}, \ldots, x^{(m)} \}$
and $S \inc \{1,\ldots,m \}$ of size $|S| = n$,
all the entries of $M_S^{-1} b(x) \in \bR^n$ are nonzero.
\elem

\bpf
Let $c \in \bR^n$ stand for $M_S^{-1} b(x)$.
The identity $M_S \, c = b(x)$ reads $\sum_{i \in S} c_i v(x^{(i)}) = v(x)$ for all $v \in \cV$.
For $j \in S$, specifying the latter when $v$ is the fundamental Lagrange interpolator $L_j$
with zeros at $x^{(i)}$,  $i \in S \setm \{j\}$, and  equal to one at $x^{(j)}$ 
yields  $c_j = L_{j}(x)$,
which is nonzero. 
\epf

The above-mentioned characterization of $n$-sparse solutions to \eqref{MinSpe}-\eqref{MinSpe2}, stated next,  is inspired by the simplex algorithm's certificate of optimality.

\bprop
\label{Prop1}
Let $\cV$ be an $n$-dimensional Chebyshev subspace of $\cC[-1,1]$.
For $x \not\in \{x^{(1)}, \ldots, x^{(m)} \}$
and $S \inc \{1,\ldots,m \}$ of size $|S| = n$,
the vector $a^{(S)}(x) \in \bR^{m}$ 
defined by $[a^{(S)}(x)]_S = M_S^{-1} b(x)$
and $[a^{(S)}(x)]_{S^c} = 0$
is a solution to \eqref{MinSpe2}
if and only if
\begin{equation}
\label{Chara}
\| M_{S^c}^{\top} M_{S}^{-\top} {\rm sgn}(M_S^{-1} b(x)) \|_\infty \le 1 .
\end{equation}
\eprop

\bpf
For simplicity of notation,
the dependence on $x$ is removed throughout the proof,
so one writes $a^{(S)}$ instead of $a^{(S)}(x)$ and $b$ instead of $b(x)$.

Suppose on the one hand that \eqref{Chara} holds. 
Then, for any $a \in \bR^m$ such that $M a  = b$,
in view of $a_S = M_S^{-1}(b - M_{S^c} a_{S^c})$, 
one has
\begin{align*}
\|a\|_1 -\| a^{(S)} \|_1
& = \|a_{S^c}\|_1 + \|a_S\|_1 - \|M_S^{-1} b\|_1
=  \|a_{S^c}\|_1 + \|M_S^{-1}(b - M_{S^c} a_{S^c})\|_1 - \|M_S^{-1} b\|_1\\
& \ge \|a_{S^c}\|_1 + \langle \sgn(M_S^{-1} b), M_S^{-1}(b - M_{S^c} a_{S^c})\rangle - \langle \sgn(M_S^{-1} b),  M_S^{-1} b \rangle \\
& = \|a_{S^c}\|_1 - \langle \sgn(M_S^{-1} b),  M_S^{-1}M_{S^c} a_{S^c} \rangle
= \|a_{S^c}\|_1 - \langle M_{S^c}^{\top} M_{S}^{-\top} (\sgn(M_S^{-1} b)),  a_{S^c} \rangle\\
& \ge \|a_{S^c}\|_1 - \| M_{S^c}^{\top} M_{S}^{-\top} (\sgn(M_S^{-1} b)) \|_\infty  \|a_{S^c} \|_1 \ge 0.
\end{align*} 
This means that $\|a\|_1 \ge \|a^{(S)}\|_1$ for any feasible vector $a \in \bR^m$ in \eqref{MinSpe2},
i.e., that $a^{(S)}$ is indeed a solution to \eqref{MinSpe2}.

Suppose on the other hand  that \eqref{Chara} does not hold.
One considers an index $\ell \in S^c$ such that $| ( M_{S^c}^{\top} M_{S}^{-\top} {\rm sgn}(M_S^{-1} b) )_\ell | > 1$.
Then, for $t \in \bR$, one defines a vector $a \in \bR^m$ satisfying $Ma = b$ via
$$
a_{S^c} = t e_\ell
\qquad \mbox{and} \qquad
a_S = M_S^{-1}(b - t M_{S^c} e_\ell).
$$ 
Since the entries of $M_S^{-1}b$ are all nonzero by Lemma \ref{LemNZ},
one has $\sgn(a_S) = \sgn(M_S^{-1}b)$ when $|t|$ is small enough,
in which case
\begin{align*}
\|a\|_1 & = \|a_S\|_1 + \|a_{S^c}\|_1
= \langle \sgn(a_S), a_S \rangle + |t|
= \langle  \sgn(M_S^{-1}b) ,    M_S^{-1}(b - t M_{S^c} e_\ell) \rangle + |t| \\
& = \|M_S^{-1} b\|_1 - t \langle M_{S^c}^{\top} M_{S}^{-\top} {\rm sgn}(M_S^{-1} b) , e_\ell \rangle  + |t| 
= \|a^{(S)} \|_1  - t  \, ( M_{S^c}^{\top} M_{S}^{-\top} {\rm sgn}(M_S^{-1} b) )_\ell  + |t|.
\end{align*} 
Thus, when $t \not= 0$ is small  enough in absolute value and chosen of the appropriate sign,
one obtains $\|a\|_1< \|a^{(S)}\|_1$, 
so that $a^{(S)}$ is not a solution to \eqref{MinSpe2}.
\epf

All the ingredients are now in place to  complete the justification of the two remaining points,
treated in reverse order of appearance.\vspace{-5mm}

\paragraph{Optimality.} 
For $k \in \{1,\ldots,m+1\}$,
the fact that the index set $S_k \inc \{1,\ldots,m\}$ of size $|S_k|=n$
is the support of a minimizer of $\|a\|_1$ subject to $M a = b(z^{(k)})$
certifies, by Proposition~\ref{Prop1},
that $\| M_{S_k^c}^{\top} M_{S_k}^{-\top} {\rm sgn}(M_{S_k}^{-1} b(z^{(k)})) \|_\infty \le 1$. 
But according to Lemma~\ref{LemNZ},  as $x$ moves through the subinterval $(x^{(k-1)},x^{(k)})$,
none of the entries of $M_{S_k}^{-1} b(x)$ can vanish,
meaning that ${\rm sgn}(M_{S_k}^{-1} b(x))$ stays the same as ${\rm sgn}(M_{S_k}^{-1} b(z^{(k)}))$,
hence implying that $\| M_{S_k^c}^{\top} M_{S_k}^{-\top} {\rm sgn}(M_{S_k}^{-1} b(x)) \|_\infty \le 1$
for any $x \in (x^{(k-1)},x^{(k)})$.
By Proposition~\ref{Prop1} again,
this ensures that
the vector supported on $S_k$ and equal to $M_{S_k}^{-1} b(x)$ there---this is precisely $a^\sharp(x)$ as defined in Algorithm \ref{Algo}---is a solution to \eqref{MinSpe}-\eqref{MinSpe2} for any $x \in (x^{(k-1)},x^{(k)})$ and any $k \in \{1,\ldots,m+1\}$.

Moreover, one also notes that $a^\sharp(x^{(k)})$---defined as $a^\sharp(x^{(k)}) = e_k$---is a solution to \eqref{MinSpe}-\eqref{MinSpe2} for $x=x^{(k)}$.
Indeed, it meets the constraint in~\eqref{MinSpe}
and its $\ell_1$-norm is $\|e_k\|_1 = 1$,
while any $a \in \bR^m$ meeting the constraint in~\eqref{MinSpe}, in particular with $v=1$, obeys $\sum_{i=1}^m a_i = 1$,
so its $\ell_1$-norm satisfies $\|a\|_1 \ge 1$.

All in all,  as outlined in Section \ref{SecOR},
the fact that the vectors $a^\sharp(x)$ are solutions to \eqref{MinSpe}-\eqref{MinSpe2} for all $x \in [-1,1]$ guarantees that $\Delta^\sharp(y) = \sum_{i=1}^m y_i a^\sharp_i$ defines a recovery map that is optimal for the full approximation problem over $\cK_{\cV,\eps}$,
provided contunuity of $a^\sharp_1,\ldots,a^\sharp_m$ can be established, which is done next.\vspace{-5mm}\qed

\paragraph{Continuity.}
Since the vector $a^\sharp(x) \in \bR^m$ was defined for $x \in (x^{(k-1)},x^{(k)})$ by $[a^\sharp(x)]_{S_k^c} = 0$ and $[a^\sharp(x)]_{S_k} = M_{S_k}^{-1} b(x)$,
the map $x \mapsto a^\sharp(x)$ is readily continuous on each subinterval $ (x^{(k-1)},x^{(k)})$.
To~ensure overall continuity,  one should check continuity at $x^{(1)},\ldots,x^{(m)}$.
Given $k \in \{1,\ldots,m\}$,
one shall verify e.g. that $a^\sharp(x)$ tends to $ a^\sharp(x^{(k)})$ as $x$ tends to $x^{(k)}$ while belonging to $(x^{(k-1)},x^{(k)})$.

From above, it is known that $a^\sharp(x^{(k)}) = e_k$ is a solution to \eqref{MinSpe}-\eqref{MinSpe2} with $x=x^{(k)}$.
One claims that the vector $\ol{a} := \lim_{x \nearrow x^{(k)}} a^\sharp(x)$ is also a solution to \eqref{MinSpe}-\eqref{MinSpe2} with $x=x^{(k)}$.
Indeed, the equality $M a^\sharp(x) = b(x)$ passes to the limit as $x \nearrow x^{(k)}$ to give $M \ol{a} = b(x^{(k)})$, so that $\ol{a} \in \bR^m$ is feasible for~\eqref{MinSpe2}.
Moreover, 
for any $a \in \bR^m$ satisfying $Ma = b(x^{(k)})$,
let $\wt{a} \in \bR^m$ be defined by $\wt{a}_{S_k^c} = a_{S_k^c}$ and $\wt{a}_{S_k} = a_{S_k} + M_{S_k}^{-1}(b(x) - b(x^{(k)}))$.
From $M_{S_k^c} \wt{a}_{S_k^c} = M_{S_k^c} a_{S_k^c}$ and
$M_{S_k} \wt{a}_{S_k} 
=  M_{S_k} a_{S_k} + b(x) - b(x^{(k)})$, 
one obtains 
$M \wt{a} = M_{S_k^c} \wt{a}_{S_k^c} + M_{S_k} \wt{a}_{S_k} = M a +  b(x) - b(x^{(k)}) = b(x)$.
By the minimality property of $a^\sharp(x)$ established before, one deduces that $\|\wt{a}\|_1 \ge \|a^{\sharp}(x)\|_1$,
and letting $x$ tend to $x^{(k)}$ yields $\|a\|_1 \ge \|\ol{a}\|_1$,
showing that $\ol{a}$ is indeed a solution to \eqref{MinSpe}-\eqref{MinSpe2} with $x=x^{(k)}$.

Since the two minimizers $e_k$ and $\ol{a}$ must have the same $\ell_1$-norm, one has
$\|\ol{a}\|_1 = 1$.
But it also holds that $\sum_{i=1}^m \ol{a}_i = 1$ because $\ol{a}$ satisfies the constraint in \eqref{MinSpe}, in particular with $v=1$.
This implies that $\ol{a}_i \ge 0$ for all $i \in \{1,\ldots,m\}$.
Besides, given $j \in S_k$,
the constraint in \eqref{MinSpe} now written when $v$ is the fundamental Lagrange interpolator  $L_j$ with zeros at $x^{(i)}$, $i \in S_k \setm \{j\}$, and equal to one at $x^{(j)}$
yields $\ol{a}_j = L_{j}(x^{(k)})$.
Thus,  if $k \not\in S_k$,
the sequence $(L_{j}(x^{(k)}))_{j \in S_k}$ would keep a constant (positive) sign.
As pointed in Section~\ref{SecChe},
this is impossible under the assumption $n \ge 3$.
It has therefore been established that $k \in S_k$---in other words,
the index of the right endpoint (and of the left one by a similar argument) of the $k$th subinterval
belongs to the support associated with this subinterval.

Finally,  the desired conclusion $\ol{a} = e_k$
follows from the fact that $\ol{a}$ and $e_k$ are now known to both be supported on $S_k$
and from the equality $M_{S_k} \ol{a} = M_{S_k} e_k$ inferred from  $M \ol{a} = b(x^{(k)}) = M e_{k} $.\vspace{-5mm} \qed

\section{Concluding remarks}
\label{SecRks}

Now that the validity of Algorithm \ref{Algo} is fully justified,
a few comments will be beneficial to put the result in perspective.\vspace{-5mm}

\paragraph{Global vs local optimality.}
The recovery map $\Delta^\sharp$ produced by Algorithm \ref{Algo} is {\em globally} optimal,
in the sense that it minimizes the {\em global} worst-case error \eqref{GWCE}.
There is also the notion of {\em local} worst-case error,
defined at a particular $y \in \bR^m$ by
$$
{\rm lwce(\Delta,y)} = \sup_{\substack{f \in \cK \\ \La_{\bm X}f = y}} \|f-\Delta(y)\|_{C(\cX)}.
$$
A {\em locally} optimal recovery map is one that assigns, to each $y \in \bR^m$,
a minimizer over $g \in C(\cX)$ of $\sup\{  \|f-g\|_{C(\cX)} , f \in \cK_y \} $ where $\cK_y := \{ f \in \cK: \La_{{\bm X}}f = y \}$,
i.e.,  a Chebyshev center of the set $\cK_y$.
This note makes no claim about local optimality.
Note that a locally optimal recovery map may involve a costly computation at each $y \in \bR^m$,
while the cost of constructing the globally optimal recovery map $\Delta^\sharp$ can be offloaded to an offline stage producing $a_1^\sharp,\ldots,a_m^\sharp$,
after which the computation of $\Delta^\sharp(y) = \sum_{i=1}^m y_i a_i^\sharp$ is almost immediate.
Note also that the recovery map $\Delta^\sharp$ is actually independent of $\eps >0$
and that it can be---at least abstractly---constructed knowing only the $v_j(x^{(i)})$ but without explicit expressions for the $v_j(x)$ themselves,
until one requires an evaluation of $\Delta^\sharp(y)$ at some point $x \in [-1,1]$.\vspace{-5mm}

\paragraph{Data and model consistency.}
The optimal recovery map $\Delta^\sharp$ put forward in this note is data-consistent---using another jargon, it is interpolatory.
Indeed, for any $y \in \bR^m$, one can see that $\Delta^\sharp(y)(x^{(\ell)}) = y_\ell$ for all $\ell \in \{1,\ldots,m\}$ from
$$
\bigg( \sum_{i=1}^m y _i a^\sharp_i \bigg)(x^{(\ell)}) = \sum_{i=1}^m y _i a^\sharp_i (x^{(\ell)})
= \sum_{i=1}^m y_i \delta_{i,\ell} = y_\ell.
$$
However, it is not model-consistent---in other words, $\Delta^\sharp(y)$ does not always belong to $\cK_{\cV,\eps}$.
Indeed,
since $\Delta^\sharp(y)$ does not depend on $\eps$,
letting $\eps \searrow 0$ in the inequality ${\rm dist}_{C(\cX)}(\Delta^\sharp(y),\cV) \le \eps$ would imply that $\Delta^\sharp(y) \in \cV$.
This is not the case, but notice that $\Delta^\sharp(y)$ is nonetheless made of pieces from $\cV$.\vspace{-5mm}

\paragraph{Streaming data.}
The cost of Algorithm \ref{Algo} is concentrated mostly on the  solutions to  about  $m$ linear programs,
which can be prohibitive for large $m$.
In the common situation of observation points $x^{(i)}$
arriving sequentially together with the values $y_i$,
it is natural to wonder whether the work done for the construction of an optimal recovery map
based on $x^{(1)},\ldots,x^{(m)}$ can be leveraged to facilitate
the construction based on $x^{(1)},\ldots,x^{(m)}$ and an added $x^{(m+1)}$---not belonging to $(x^{(m)},1)$ in the following discussion.
When creating the new $S_k$---to fix the ideas, 
a support $S$ associated with a subinterval $I$---a simple idea is to provide the simplex algorithm with a `warm start',  i.e.,  a good initial guess $\wt{S}$.
For instance, if $I$ does not contain the added $x^{(m+1)}$ as an endpoint, 
then the guessed $\wt{S}$ can be chosen as the old $S$ associated with~$I$.
But if $I$ does contain $x^{(m+1)}$ as an endpoint,
e.g.  $I$ is the left part of an old subinterval $J$ split by $x^{(m+1)}$,
since $S$ should contain (the indices of) the left endpoint of $J$ and of $x^{(m+1)}$,
then the guessed $\wt{S}$ can be a slight modification of the old $S$ associated with $J$
obtained 
by removing (the index of) the right endpoint of $J$
and replacing it by  (the index of) $x^{(m+1)}$.
Empirically, 
the speed-up is modest and becomes more significant when the size $n$ of the supports gets closer to $m$. \vspace{-5mm}

\paragraph{Ratio of norms.}
The arguments underpinning Algorithm \ref{Algo} allow one to compute the exact value of the compatibity indicator $\mu_{\cV,\Id}(\La_{\bm X})$ and, as an interesting side product,  they transform into a practical numerical recipe to compute the exact value of the maximal ratio of uniform and discrete norms in~$\cV$,
i.e., of 
\be
\label{Ratio}
\rho_{\cV, {\bm X}} := \max_{v \in \cV \setm \{ 0\}} \f{\|v\|_\infty}{\max_{i \in \{1,\ldots,m\}} |v(x^{(i)})|}.
\ee
When $\cV = \cP_n$ is the space of algebraic polynomials of degree $<n$,
this maximal ratio has been well studied---in particular, for equispaced points, see \cite{CopRiv,Rak}---and a computational method akin to Remez algorithm
 has been proposed in \cite[Section 6]{APS}.
The recipe uncovered here for a Chebyshev space $\cV$ of dimension $n \ge 3$ and containing the constant functions amounts to solving the $m+1$ linear programs from Algorithm \ref{Algo} to create the index sets $S_1,\ldots,S_{m+1}$.
After that, it is summarized as
\be
\label{AlgoMu}
\rho_{\cV, {\bm X}} = \max_{k \in \{1,\ldots,m+1 \} }  \max_{x \in [x^{(k-1)},x^{(k)}]}
 \sum_{j \in S_k} {\rm sgn}( ( M_{S_k}^{-1}  b(z^{(k)}))_j ) ( M_{S_k}^{-1}  b(x))_j . 
\ee
In short, one needs to compute the maximal values of functions from $\cV$ on $m+1$ subintervals.
The justification of \eqref{AlgoMu} relies on the identity $\mu_{\cV,\Id}(\La_{\bm X}) = 1 + \rho_{\cV,{\bm X}}$ established in \cite[Lemma 6.2]{DFPW} and on the following observation:
\begin{align*}
\mu_{\cV,\Id}(\La_{\bm X}) & = \sup_{h \in \ker \La_{\bm X}  \setm \{0 \}} \f{\|h\|_{C[-1,1]}}{{\rm dist}_{C[-1,1]}(h,\cV)}
= \max_{k \in \{1,\ldots,m+1\} } \max_{x \in [x^{(k-1)},x^{(k)}]}
\sup_{h \in \ker \La_{\bm X}  \setm \{0 \}} \f{|h(x)|}{{\rm dist}_{C[-1,1]}(h,\cV)}\\
& =  \max_{k \in \{1,\ldots,m+1\} } \max_{x \in [x^{(k-1)},x^{(k)}]} \big( \mu_{\cV,\delta_x}(\La_{\bm X}) \big)
= \max_{k \in \{1,\ldots,m+1\} } \max_{x \in [x^{(k-1)},x^{(k)}]} \big( 1 + \|a^\sharp(x)\|_1 \big) \\
& = 1+ \max_{k \in \{1,\ldots,m+1\} } \max_{x \in [x^{(k-1)},x^{(k)}]} \sum_{j \in S_k}
 | ( M_{S_k}^{-1}  b(x))_j |\\
& = 1+  \max_{k \in \{1,\ldots,m+1\} } \max_{x \in [x^{(k-1)},x^{(k)}]} \sum_{j \in S_k}
 {\rm sgn}( ( M_{S_k}^{-1}  b(z^{(k)}))_j ) ( M_{S_k}^{-1}  b(x))_j ,
\end{align*}
where the last step made use of the fact that the signs of $M_{S_k}^{-1}  b(x)$ do not change throughout the interval $[x^{(k-1)},x^{(k)}]$.
Of course, the practicality of the recipe stemming from \eqref{AlgoMu} depends on the ability to compute the maxima of functions from $\cV$.
This task can be efficiently performed in Chebfun \cite{chebfun},
an open-source {\sc matlab} package for numerical computations with functions.
This feature,  together with the easy handling of piecewise functions, 
explains why Chebfun was preferred for the implementation of Algorithm \ref{Algo} in the reproducible file accompanying this note (available on the author's webpage).
Furthermore, if $\cV$ is a space of trigonometric or algebraic polynomials,
then the maximum on a subinterval of a function from $\cV$ can be computed via semidefinite programming, as explained in \cite{FouPow}, see Theorem 3.1 in particular. \vspace{-5mm}

\rev{
\paragraph{Best choice of evaluation points.}
Throughout this note,
the points $x^{(1)},\ldots,x^{(m)}$ were prescribed.
If one could select them freely,
one would naturally want to do so in a way that makes the minimal worst-case error as small as possible.
According to the previous considerations,
this amounts to minimizing over all pointsets ${\bm X}_m$ of size $m$ the ratio $\rho_{\cV, {\bm X}_m}$ introduced in \eqref{Ratio}.
This is bound to be a difficult problem,
as it is unresolved even for $m=n$ and $\cV = \cP_n$.
Indeed, it is not hard to see that  $\rho_{\cP_n, {\bm X}_n}$ coincides with the $L_\infty$-operator norm of the interpolation operator at ${\bm X}_n = \{x^{(1)},\ldots,x^{(n)}\}$,
aka the Lebesgue constant.
Pointsets with nearly optimal Lebesgue constant are known explicitly,
but  pointsets with genuinely optimal Lebesgue constant are not,
even though they have been characterized a long time ago, see \cite{deBPin,Kil}. 
\vspace{-5mm}
}

\rev{
\paragraph{Other observation functionals.}
The $\ell_1$-minimizations at the heart of Algorithm \ref{Algo} appear thanks to the presence of point evaluations.
The situation becomes more complicated if arbitrary observation functionals were involved.
Nonetheless, if one had the freedom to use any observation functionals,
it would be natural to wonder about the power of point evaluations:
is the minimal worst-case error much smaller with unrestricted functionals than with point evaluations only?
No attempts were made to answer this question in the context of this note,
but some recent advances are worth pointing out in a related context where the recovery performance is assessed via the $L_2$-norm rather than the uniform norm in \eqref{GWCE}, see e.g.  \cite{DKU} and the references therein.
\vspace{-5mm}
}


\end{document}